\documentclass{conm-p-l}
\usepackage{amssymb}
\usepackage{epsf}

\usepackage{colordvi}
\usepackage{graphicx}

\numberwithin{equation}{section}
\newtheorem{theorem}{Theorem}[section]
\newtheorem{defn}[theorem]{Definition}

\newtheorem{example}[theorem]{Example}
\newtheorem{remark}[theorem]{Remark}

\theoremstyle{remark}

\newcounter{FNC}[page]
\def\fauxfootnote#1{{\addtocounter{FNC}{2}\Magenta{$^\fnsymbol{FNC}$}%
     \let\thefootnote\relax\footnotetext{\Magenta{$^\fnsymbol{FNC}$#1}}}}

\def\cocoa{{\hbox{\rm C\kern-.13em o\kern-.07em C\kern-.13em o\kern-.15em A}}}

\newcommand{\C}{{\mathbb{C}}}

\renewcommand{\P}{{\mathbb{P}}}
\newcommand{\R}{{\mathbb{R}}}

\newcommand{\DeCo}[1]{\Blue{#1}}

\title[Experimentation at the Frontiers of Reality in Schubert Calculus]{Experimentation
  at the Frontiers of Reality\\ in Schubert Calculus} 

\author[Hillar]{Christopher Hillar}
\address{Christopher Hillar\\The Mathematical Sciences Research Institute\\
         17 Gauss Way\\
         Berkeley, CA 94720-5070}
\email{chillar@msri.org}
\author[Garc{\'\i}a]{Luis Garc\'ia-Puente}
\address{Luis Garc\'ia-Puente\\Department of Mathematics and Statistics\\
         Sam Houston State University\\
         Huntsville\\
         TX \ 77341}
\email{lgarcia@shsu.edu}
\urladdr{www.shsu.edu/\~{}ldg005}
\author[Mart\'in del Campo]{Abraham Mart\'in del Campo}
\address{Abraham Mart\'in del Campo\\Department of Mathematics\\
         Texas A\&M University\\
         College Station\\
         TX \ 77843}
\email{asanchez@math.tamu.edu}
\urladdr{http://www.math.tamu.edu/\~{}asanchez}
\author[Ruffo]{James Ruffo}
\address{James Ruffo\\Department of Mathematics, Computer Science, \& Statistics\\
         State University of New York\\
         College at Oneonta\\
         Oneonta, NY 13820}
\email{ruffojv@oneonta.edu}
\urladdr{http://employees.oneonta.edu/ruffojv/}
\author[Teitler]{Zach Teitler}
\address{Zach Teitler\\Department of Mathematics\\
         Texas A\&M University\\
         College Station\\
         TX\ 77843}
\email{zteitler@math.tamu.edu}
\urladdr{http://www.math.tamu.edu/\~{}zteitler}
\author[Johnson]{Stephen L.~Johnson}
\address{Stephen L.~Johnson\\Department of Mathematics\\
         Texas A\&M University\\
         College Station\\
         TX \ 77843}
\email{steve@math.tamu.edu}
\urladdr{http://www.math.tamu.edu/\~{}steve.johnson}
\author[Sottile]{Frank Sottile}
\address{Frank Sottile\\Department of Mathematics\\
         Texas A\&M University\\
         College Station\\
         TX \ 77843}
\email{sottile@math.tamu.edu}
\urladdr{http://www.math.tamu.edu/\~{}sottile}

\thanks{Research of Sottile supported in part by NSF grant DMS-070105}
\thanks{Research of Hillar  supported in part by an NSF Postdoctoral Fellowship and an NSA
  Young Investigator grant} 

\begin{document}

\begin{abstract}
 We describe a general framework for large-scale computational experiments in mathematics
 using computer resources that are available in most mathematics departments.
 This framework  was developed for an experiment 
 that is helping to formulate and test conjectures in the real Schubert calculus. 
 Largely using machines in instructional computer labs during off-hours and University breaks,
 it consumed in excess of 350 GigaHertz-years of computing in its first six months of
 operation, solving over 1.1~billion polynomial systems.
\end{abstract}

\maketitle
%
\section*{Introduction}
Mathematical discovery has long been informed by experimentation and
computation. 
Understanding key examples is typically the first step towards formulating theorems and
devising proofs. 
The computer age enables many more potentially intricate examples to be studied
than ever before.
Sometimes, this leads to a fruitful dialog between theory and experiment.
Other times, this work leads serendipitously to new ideas and
theorems.  
Many examples are described in the books~\cite{BBCGLM,BB08}.

We believe there is much greater potential for computer-aided experimentation than what
has been achieved.
This is particularly true for scientific discovery, using advanced computing to study 
subtle phenomena and amass evidence for the
mathematical facts which will become the theorems of tomorrow.
Currently, much computer experimentation is (often appropriately) on a fairly small scale. 
A notable exception 
is Odlyzko's study~\cite{Od87} (using Cray supercomputers) of the zeroes
of Riemann's $\zeta$-function on the critical line $\frac{1}{2}+\R\sqrt{-1}$, which
led to a rich data set that has
stimulated much intriguing mathematics~\cite{Diaconis}.

A different large scale use of computers is the Great Internet Mersenne Prime
Search (GIMPS)~\cite{GIMPS}, which searches for primes of the form $2^p{-}1$
for $p$ a prime, such as $3,7,31$, and $127$.
Volunteers run software on otherwise
idle computers to search for Mersenne primes.
This project has found the largest known primes since it started in 1996.
Daily, it uses over 60 GigaHertz-years of computation.

GIMPS is a mathematical analog of big-science physics.
We feel there is more scope for such investigations in mathematics.
We describe our use of a supercomputer to
study a conjecture in the real Schubert calculus,
which may serve as a model for research in mathematics based on computational experiments.
Rather than Odlyzko's Cray supercomputers, or GIMPS's thousands of volunteers,
we use more pedestrian computer 
resources that are available to many mathematics departments together with modern 
(and free) software tools such as Perl~\cite{Perl}, MySQL~\cite{mysql}, and PHP~\cite{php},
as well as freely available mathematical software such as
Singular~\cite{Singular}, Macaulay 2~\cite{M2}, and Sage~\cite{sage}.

This is a methods paper whose purpose is to explain the framework we developed.
We do not present mathematical conclusions from
this ongoing computational experiment, but instead explain how you, the reader, can take
advantage of readily available yet often underutilized computer resources to employ in
your mathematical research. 

To get an idea of the available resources, in its first six
months of data acquisition, this experiment used over 350 GigaHertz-years
of computing primarily on 191 computers in instructional labs that are maintained
by the Department of Mathematics at Texas A\&M University~\cite{Calclabs}.
When the labs are not in use, the machines become a cluster computing resource that
provides over 500 computational cores for a peak performance of 1.971 Teraflops with 296GB of
total memory.
This experiment uses a supercomputer
moonlighting from its day job of calculus instruction.

The authors of this note include Johnson, who configured the labs as a Beowulf cluster,
enabling their use for this computation.
Our software 
was written and maintained by the remaining authors,
who include current and former postdocs and students working with Sottile.
We are organized into a vertically-integrated team where 
the senior members work with and mentor the junior members.
The overall software design and much of its
implementation is due to Hillar.

A key feature of this experiment is its robustness---it can and has recovered from many
faults, including emergency system shutdown, database erasure, inexplicable computer
malfunction, as well as day-to-day network failures. 
It is also repeatable, using a pseudorandom number generator with fixed seeds.
This repeatability will allow us 
to rerun a large segment of our calculations on a different
supercomputer~\cite{Brazos} using different mathematical software than the initial run.
This will be an unprecedented test of the efficacy of different implementations of 
our basic mathematical algorithms of Gr\"obner basis computation and real root counting.

This experiment is part of a long-term study of a striking conjecture
in the real Schubert calculus made by Boris Shapiro and Michael Shapiro in 1993.
This includes two previous large computational experiments~\cite{RSSS,So_Shap}
(and several smaller ones~\cite{RoSo,So_fulton}), as well as more traditional
work~\cite{EG_02,EGSV,KhSo,So_special,So_quantum}, including proofs of the Shapiro Conjecture
for Grassmannians~\cite{MTV_Sh,MTV_R}.
This story was the subject of a Current Events Bulletin Lecture at
the 2009 AMS meeting and a forthcoming article in the AMS Bulletin~\cite{CEB}.

This experiment is possible only because we may model the geometric problems we study 
on a computer, efficiently solve them, record the results, and automate this process.
We describe some background in Section~\ref{Sec:Shapiro} and
the mathematics of the  computations in Section~\ref {S:math}.
In Sections~\ref{S:Resources}--\ref{S:quality}, we explain the resources (human, hardware,
and software) we utilized, the architecture  of the experiment, how we ran it on a
cluster, and the measures that we took to maintain the quality of our data.
We end with some conclusions and remarks.


%
\section{The Shapiro Conjecture and Beyond}\label{Sec:Shapiro}

Our goal is to describe the design and execution of a large scale computation,
which may serve as a model for other experiments in mathematics.
While many aspects of our experiment are universal, the details are specific to the
questions we are studying.
We give some mathematical background to provide context.

Some solutions to a system of real polynomial equations are real and the rest
occur in complex conjugate pairs.
While the structure of the equations determines the total number of solutions, 
the distribution between the two types depends subtly on the coefficients.
Surprisingly, sometimes there is additional structure which leads to
finer information in terms of upper bounds~\cite{BBS,Kh91} or lower bounds~\cite{EG01,SS}
on the number of real solutions.
The Shapiro Conjecture is the extreme situation of having {\sl only} real solutions.

\begin{figure}[htb]
 \[
  \begin{picture}(310,160)(0,4)
   \put(0,0){\includegraphics[height=175pt]{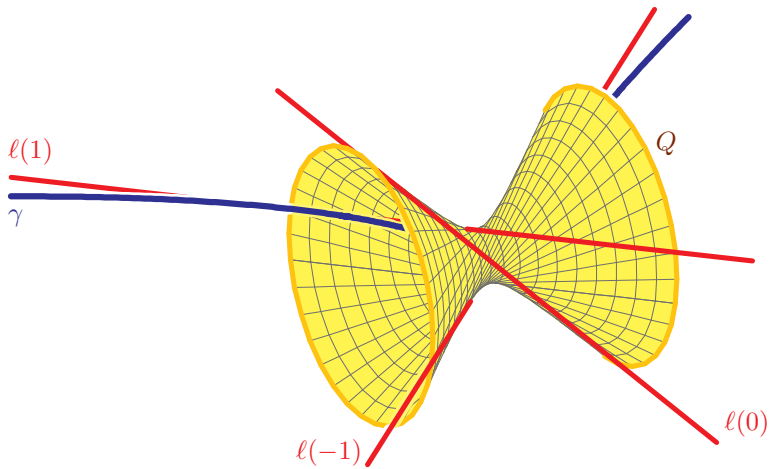}}
   \put(112,3){\Red{$\ell(-1)$}} \put(274,14){\Red{$\ell(0)$}}
   \put(3,117){\Red{$\ell(1)$}} \put(3,93){\Blue{$\gamma$}}
   \put(248,120){\Brown{$Q$}}
  \end{picture}\vspace{-10pt}
 \]
 \caption{Quadric containing three lines tangent to the curve $\gamma$\label{F:TanQuad}.} 
 \end{figure}

We give an example.
Set $\gamma(t):=(6t^2{-}1, \frac{7}{2}t^3{+}\frac{3}{2}t, -\frac{1}{2}t^3{+}\frac{3}{2}t)$,
which is a curve $\DeCo{\gamma}\colon\R\to\R^3$.
We ask for the finitely many lines that meet four tangent lines to $\gamma$, which we take
to be tangent at the points $\gamma(t)$ for $t=-1,0,1$, and some point 
$s\in(0,1)$.  
The first three tangents lie on the quadric $Q$ defined by  $x^2-y^2+z^2=1$.
We show this in Figure~\ref{F:TanQuad}, where $\ell(t)$ is the tangent line
at the point $\gamma(t)$. 
 \vspace{-5pt}

\begin{figure}[htb]
 \[
  \begin{picture}(254, 116)(-14,4)
   \put(0,0){\includegraphics[height=110pt]{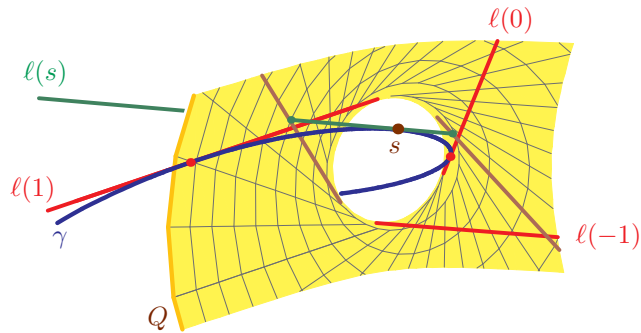} }
   \put(-14,51){\Red{$\ell(1)$}} \put(200,33){\Red{$\ell(-1)$}} 
   \put(167,115){\Red{$\ell(0)$}}
   \put(129,67){\Brown{$s$}}
   \put(38,2){\Brown{$Q$}}
   \put(2,33){\DeCo{$\gamma$}} \put(-9,95){\ForestGreen{$\ell(s)$}}
  \end{picture}\vspace{-10pt}
\] 
\caption{Configuration in the throat of the quadric\label{F:throat}.}
 \end{figure}

These first three tangents lie on one ruling of $Q$ and the lines in the other ruling 
are those meeting them.
Lines meeting all four tangents correspond to the (two) {\it a priori} complex points
where the fourth tangent $\ell(s)$ meets the quadric.
As we see in Figure~\ref{F:throat}, for any $s\in(0,1)$, $\ell(s)$ meets the quadric in two
real points, giving two real solutions to this instance of the problem of four lines.

The Schubert calculus~\cite{Fu97,FuPr} asks for the 
linear spaces that have specified positions with respect to
other, fixed (flags of) linear spaces.
For example, what are the 3-planes in $\C^7$ meeting 12 given 4-planes non-trivially? 
(There are 462~\cite{Sch1886c}.)
The specified positions are a \DeCo{{\sl Schubert problem}}, for example, the
Schubert problem of lines meeting four lines in 3-space.
The fixed linear spaces imposing the conditions give 
an \DeCo{{\sl instance}} of the Schubert problem, so that the lines 
$\ell(-1)$, $\ell(0)$, $\ell(1)$, and $\ell(s)$ give an instance of the problem of 
four lines.
The number of solutions depends upon the Schubert problem, while the solutions
depend upon the instance.

The Shapiro Conjecture begins with a rational normal curve $\gamma\colon\R\to\R^n$, which is
any curve projectively equivalent to the moment curve,
\[
    t\ \longmapsto\ (t,\, t^2,\, \dotsc,\, t^n)\,.
\]
In 1993, Boris Shapiro and Michael Shapiro conjectured that if the fixed linear spaces
osculate a rational normal curve, then all solutions to the Schubert problem are real.
Initially, the statement seemed too strong to be possibly true.
This perception changed dramatically after a few computations~\cite{RoSo,So_PSPM}, 
leading to a systematic study of the conjecture for Grassmannians, both theoretical and 
experimental~\cite{So_Shap} in which about 40,000 instances were computed of 11 different
Schubert problems.
%
%
%
Several extremely large instances were also verified by others~\cite{FRZ,Ve}.

This early study led to a proof of the Shapiro Conjecture in a limiting sense for
Grassmannians~\cite{So_special} and a related result in the quantum cohomology of
Grassmannians~\cite{So_quantum}, which drew others to the area.
Eremenko and Gabrielov~\cite{EG_02} proved it for Grassmannians of codimension 2
subspaces where the Shapiro Conjecture becomes the statement
that a univariate rational function with only real critical points is (equivalent to) a
quotient of real polynomials.

Later, Mukhin, Tarasov, and Varchenko~\cite{MTV_Sh} used ideas from integrable systems 
to prove the Shapiro Conjecture for Grassmannians.
They later gave a second proof~\cite{MTV_XXX} that revealed deep
connections between geometry and representation theory.
This story was the subject of a Current Events Bulletin Lecture at the January 2009
AMS meeting and a forthcoming article in the AMS Bulletin~\cite{CEB}.

%
\subsection{Beyond the Grassmannian}

The Shapiro Conjecture makes sense for any flag manifold (compact rational homogeneous
space).
Early calculations~\cite{So_fulton} supported it for orthogonal
Grassmannians but found counterexamples for general $SL_n$-flag manifolds and
Lagrangian Grassmannians.
Calculations suggested modifications in these last two cases~\cite{So_ERAG} and limiting
versions were proved~\cite{So_fulton},
and the conjecture for the orthogonal Grassmannian wasjust proven by Purbhoo~\cite{Pur}.

The modification for $SL_n$-flag manifolds, the \DeCo{{\sl Monotone Conjecture}}, was
refined and tested in a computational experiment involving Ruffo and Sottile~\cite{RSSS}. 
That ran on computers at the University of Massachusetts, the
Mathematical Sciences Research Institute, and Texas A\&M University, using 15.76
GigaHertz-years of computing to study over 520 million instances of 
1126 different Schubert problems on 29 flag manifolds.
Over 165 million instances of the Monotone Conjecture were verified, and the
investigation discovered many new and interesting phenomena. 

For flags consisting of a codimension 2 plane lying on a hyperplane, the Monotone Conjecture 
is a special case of a statement about real rational functions which Eremenko,
et.~al~\cite{EGSV} proved. 
Their work leads to a new conjecture for Grassmannians.
A flag is \DeCo{{\sl secant along an interval $I$}} of a
curve if every subspace in the flag is spanned by its intersections with $I$.
The \DeCo{{\sl Secant Conjecture}} asserts that if the flags in a Schubert problem on a 
Grassmannian are \DeCo{{\sl disjoint}} in that they are secant along disjoint intervals of a
rational normal curve, then every solution is real. 
It is true for Grassmannians of codimension 2 subspaces, by the result of Eremenko, et.~al.

\begin{figure}[htb]
\[
  \begin{picture}(320,135)(0,-1)
   \put(0,7){\includegraphics[height=135pt,viewport=5 77 430 282,clip]{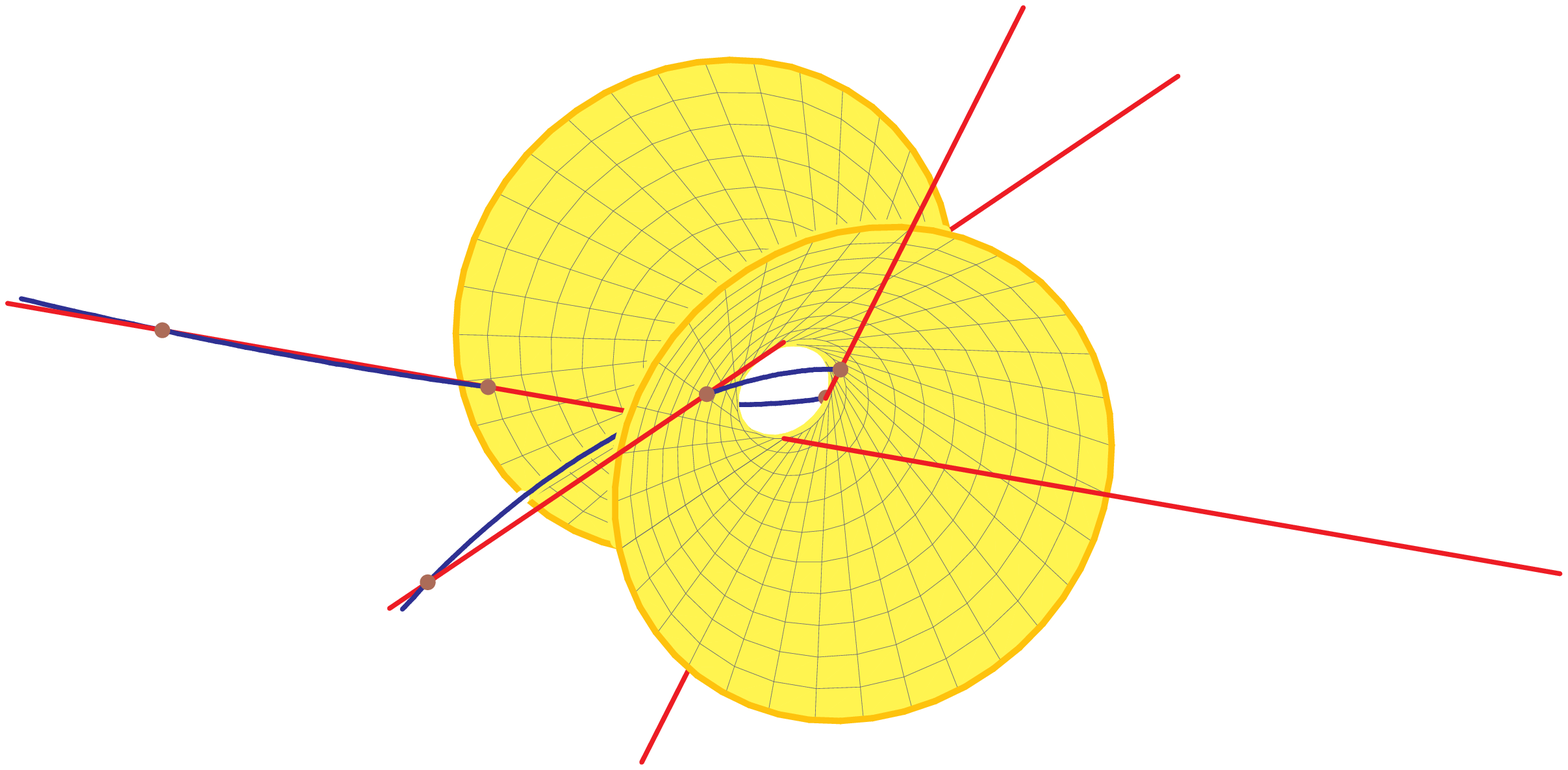}}
   \thicklines
   \put(105,50){$\gamma$}\put(115,51){\vector(2,-1){25}}
   \put(227,-2){\vector(0,1){77.8}}
  \end{picture}\vspace{-10pt}
\]
\caption{The problem of four secant lines.}
\label{F:secant}
\end{figure}

Consider this for the problem of four lines.
The hyperboloid in Figure~\ref{F:secant} contains three lines that are secant to 
$\gamma$ along disjoint intervals.
Any line secant along the indicated arc (which is disjoint from the other three intervals)
meets the hyperboloid in two points, giving two real solutions to this
instance of the Secant Conjecture.

We are testing the Secant Conjecture for many Schubert problems on small Grassmannians of
$k$-planes in $n$-space.
(See Table~\ref{T:number_problems}.)
\begin{table}[htb]
 \begin{tabular}{|c||c|c|c|c|}\hline
  $k\backslash n{-}k$ & 2&3&4&5\\\hline\hline
     2 & 1 & 5 & 22 & 81 \\\hline
     3 & 5 &63 & 94 & \\\hline
     4 &22 &  & &\\\hline
     5 & 81 & & & \\\hline
 \end{tabular}\vspace{4pt}
\caption{Schubert problems studied on $G(k,n)$ as of 20 May 2009.}
\label{T:number_problems}
\end{table}\vspace{-15pt}
Instances of the Secant Conjecture are substantially
more difficult to compute than instances of the Shapiro Conjecture.
Consequently, this experiment has used much more computing than 
experiments for the Shapiro and Monotone Conjectures.

%
\section{Solving Schubert problems}\label{S:math}
Our core mathematical task is the following:
Given a Schubert problem on a Grassmannian and secant flags to a rational normal
curve, formulate the Schubert problem as a
system of equations and determine its number of real solutions.

A Schubert problem is a list $(w_1,\dotsc,w_s)$ of 
conditions
to be imposed on $k$-planes in $n$-space.
A fixed flag instantiating the condition $w_i$ is secant along some $m_i$ points of the
rational normal curve $\gamma$, so we need $m:=m_1+\dotsb+m_s$ points of $\gamma$
for the given Schubert problem.
For the problems we study, $8\leq m\leq 52$.
Given a set $T$ of $m$ numbers, we use the points 
$\{\gamma(t)\mid t\in T\}$ to construct secant flags.

Given these secant flags, we formulate the Schubert problem in local coordinates for the
Grassmannian (see~\cite{Fu97,RSSS,So_Shap} for details), obtaining a set of equations whose
common zeroes represent the solutions to the Schubert problem in the local coordinates.
We then eliminate all but one variable from the equations, obtaining an \DeCo{{\sl eliminant}}.
The Shape Lemma~\cite{So_M2} implies that number of real roots of this eliminant equals the number of
real solutions to the original Schubert problem,
provided that it is square-free and has degree equal to the expected number of complex solutions,
which is easily checked.

This procedure counts the real solutions to a Schubert problem involving secant flags.
To compute an instance of the Secant Conjecture, we select $m$ real numbers $T$ and use them
to form disjoint secant flags.
We also compute cases when the flags are not disjoint.
To each set $T$, we compute five problems involving secant flags.
The first is an instance of the Secant Conjecture, but for each of the other four, we randomly
alter the assignment of points to get secant flags that are not necessarily disjoint.
The \DeCo{{\sl overlap number}} measures how far the flags are from disjoint.
It is zero if and only if the flags are disjoint.

For each Schubert problem, we perform these steps thousands to millions of times.
We repeatedly select subsets $T$ of $m$ numbers from a fixed set of 111 rational numbers, 
those $p/q$ with $p^2+q^2\leq 121$.
These have small arithmetic height (which affects computation speed).
\begin{figure}[htb]
\[
  \begin{picture}(150,158)(0,5)
   \put(0,7){\includegraphics[height=145pt]{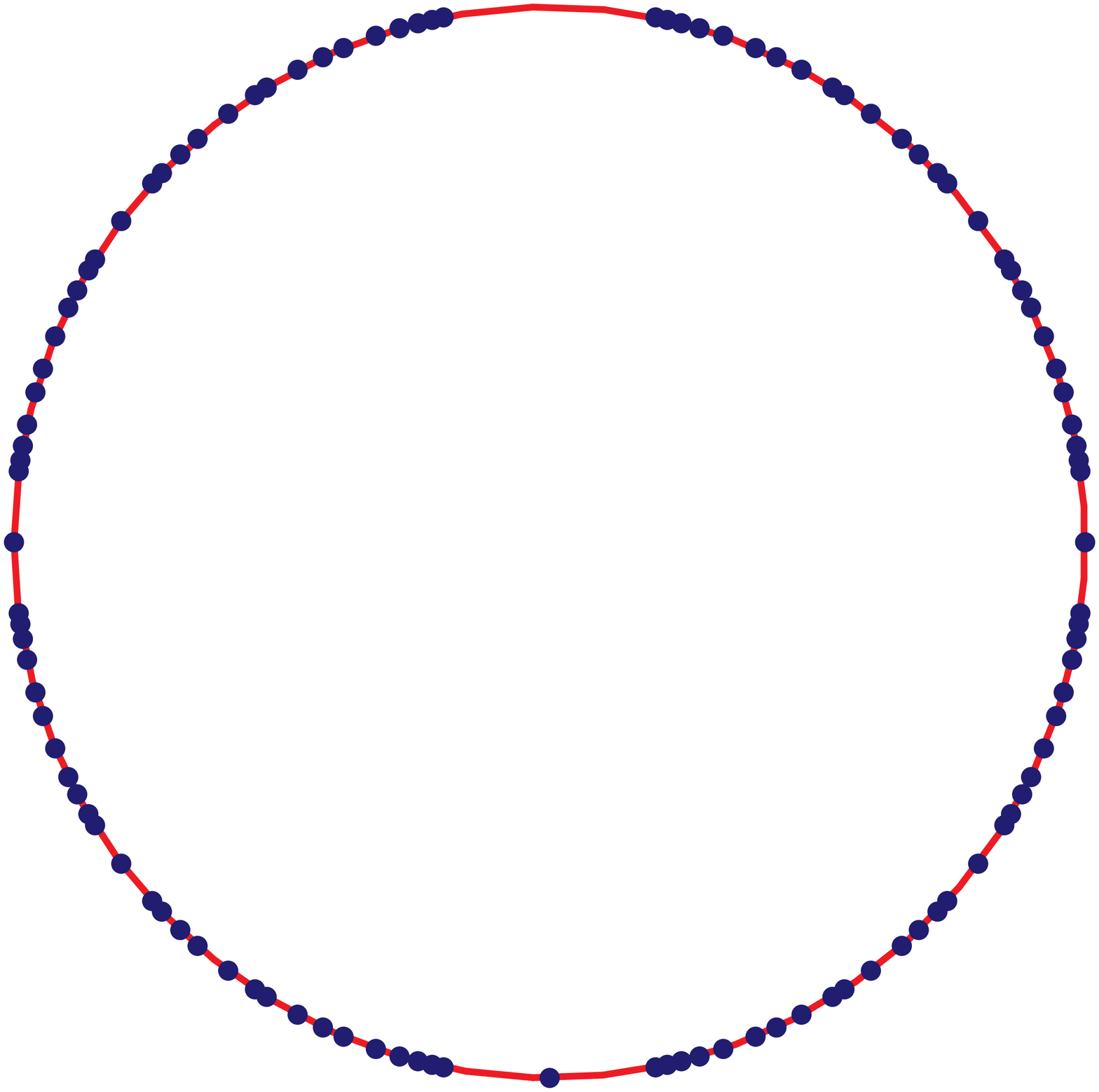}}
   \put(70.5,-2){$0$}  \put(67,154){$\infty$}
  \end{picture}
\] 
\caption{111 points along $\R\P^1$.}\label{F:RP1}
\end{figure}
Figure~\ref{F:RP1} shows these points along $\R\P^1$.

Table~\ref{T:Sample_table} shows the results
for a Schubert problem  with 16 solutions on $G(3,6)$, the Grassmannian of 3-planes in 6-space.  
There were $20,000,000$ computed instances of this problem which used 4.473 GigaHertz-years.
The rows are labeled with the even integers from 0 to 16 as the number of real solutions has
the same parity as the number of complex solutions.
The columns are labeled with overlap number, but only the first few and the summary column
are shown. 
\begin{table}[htb]
%

\noindent{\small 

\noindent\begin{tabular}{|r||r|r|r|r|r|r|r|c||r|}\hline
\textbackslash & 0&1&2&3&4&5&6&\ $\dotsb$\ &Total\\\hline\hline
0 &&&&&&&20&\ $\dotsb$\ & 7977\\\hline
2 &&&&&&&116&\ $\dotsb$\ & 88578\\\hline
4 &&&&6154&23561&526&3011&\ $\dotsb$\ & 542521\\\hline
6 &&&&25526&63265&2040&9460&\ $\dotsb$\ & 1571582\\\hline
8 &&&&33736&78559&2995&13650&\ $\dotsb$\ & 2834459\\\hline
10 &&&&25953&39252&2540&11179&\ $\dotsb$\ & 3351159\\\hline
12 &&&&35578&44840&3271&14160&\ $\dotsb$\ & 2944091\\\hline
14 &&&&17367&17180&1705&7821&\ $\dotsb$\ & 1602251\\\hline
16 &4568553&&182668&583007&468506&36983&83169&\ $\dotsb$\ &7057382 \\\hline\hline
Total &4568553&&182668&727321&735163&50060&142586&\ $\dotsb$\ &20000000 \\\hline
\end{tabular}

}\vspace{4pt}
\caption{Real solutions v.s. overlap number.}\vspace{-15pt}
\label{T:Sample_table}
\end{table}
The column with overlap number 0 represents tests of the Secant Conjecture.
Since its only entries are in the row for 16 real solutions, the Secant
Conjecture was verified in $4,568,553$ instances.
The column labeled 1 is empty because flags for this
problem cannot have overlap number 1.
The most interesting feature is that for overlap number 2, all
solutions were real, while for overlap numbers 3, 4, and 5, at least 4 of the 16
solutions were real, and only with overlap number 6 and greater does the Schubert problem have
no real solutions. 
This inner border, which indicates that the reality of the Schubrt problem does not completely
fail when there is small overlap, is found on many of the other problems that we investigated
and is a new phenomenon that we do not understand.

%
%
\section{Resources for this experiment}\label{S:Resources}

Creating the software and managing this computation is a large project.
To accomplish it, we formed a vertically-integrated team of graduate students and postdoctoral
fellows  under the direction of a faculty member and used modern software tools (Perl, MySQL,
PHP) to automate the computation as well as store and visualize data.
This software runs on many different computers, but primarily on a supercomputing
cluster whose day job is calculus instruction.

%
\subsection{People}
The authors of this note include Johnson, who created and maintains the supercomputer
we use, as well as a research team of current and former graduate students and postdoctoral fellows
who have worked with Sottile.
Pooling our knowledge in mathematics and in software development, we
shared the work of creating and running this project.
This structure enabled the senior members to mentor and train the junior members.

We modeled our team structure on the working environment in a laboratory.
This led to a division of labor and to other collaborations.
For example, Garc\'ia-Puente and Ruffo wrote the mathematical heart of the computation in a
Singular library.
Hillar, who had experience in the software industry, provided the conceptual framework
and contributed most of the Perl code.
He worked on some of this with Mart\'in del Campo, who now maintains the PHP webpages we use
to monitor the computation.
Sottile and Teitler maintain the software and the shell scripts for controlling the computation
and ensuring the integrity of the data, and Teitler rewrote the library of 
our main mathematical routines in Macaulay 2.
This project has led to unrelated research collaborations between Hillar
and Mart\'in del Campo and between Sottile and Teitler.

Our team includes two more junior members who have not yet contributed code and will soon
include an additional postdoc.
The web of collaboration and mentoring is designed to help integrate them into 
future projects.

\subsection{Hardware}\label{S:hardware}

All mathematics departments have significant, yet deeply underutilized computing resources
available in office desktop computers.
There are sociological problems that can arise, for example, when
your colleague has email problems while his computer is running your software.
While these can be overcome, there are often simpler alternatives.
Many institutions have some cluster computing resources, and there are regional and national
supercomputers available for research use~\cite{teragrid}.
Computers in instructional labs are another resource.
With sufficient interest and a modest expenditure, these can be used for research.

The Texas A\&M University mathematics department maintains computers for undergraduate
instruction.
Johnson, the departmental systems administrator, installed job scheduling software enabling their use
as a computing cluster outside of teaching hours.
The availability of this resource, as much as our mathematical interests, was the catalyst for
this computational experiment.
It has been the source of $95\%$ of the computing for this experiment,
%
%
%
which also used some desktop computers at Texas A\&M University and at Sam Houston State
University, as well as personal laptops and clusters at the homes of Garc\'ia-Puente and Sottile.

%
\subsection{Software}

The computer programming community has developed a vast library of free, open-source software
that mathematicians can use for research purposes.
The three software tools that we use the most,
other than specialized mathematical software, are Perl, MySQL, and PHP.
We selected them because of our familiarity with them and their widely available
documentation.
In addition to excellent manuals~\cite{ORPerl,ORMySQL,ORPHP},
there are many web pages showing documentation, tutorials, answers to frequently
asked questions, and pieces of code.

The distributed nature of our computation, its size, and the amount of data we store, led us to
organize the computation around a database to store the results and status of the computation.
For this, we chose MySQL, a freely available high-quality database program.
The actual database is located on a Texas A\&M University mathematics department server and may
be accessed from anywhere in the world.
In particular, we can (and do) monitor and manage the computation remotely.


Perl is a general-purpose programming language with especially strong facilities for text
manipulation and communication with other programs, including MySQL.
We use Perl to connect together the mathematical programs that actually perform our
calculations (Singular, Maple) with the database.

These data are viewed through web pages, which are dynamically generated using PHP, a
programming language designed exactly for this purpose.
Our interface for monitoring the experiment is at our project's
web page~\cite{secant}.

This model---computations on individual computers controlled by a central database---scales
well and is very flexible.
It can run on a single computer using a local database (e.g.~a
personal laptop), on a cluster at one's home or department, or on machines at different
institutions.


%
%
\section{Architecture of Computation}\label{S:Architecture}

We wanted to conduct a large computational experiment using distributed, heterogeneous computer
resources and have the computation be largely automated as well as robust, repeatable, and
reliable. 
To accomplish this, we organized it around a database that records all aspects of
the computation.
In Section~\ref{S:Calclabs} we explain how we run this on a cluster and in
Section~\ref{S:quality} we discuss measures to enhance the quality of our data.
Here, we focus on the organization of the computation in our software:
The basic mathematical procedures, interaction with the database, 
dividing this large computation into reasonable-sized pieces, and lastly 
selecting problems to compute and setting parameters of the computation.

%
\subsection{Implementation}

Our computation is split between three subsystems,
a controlling Perl script and mathematical computations in Singular and in Maple.
We explain these choices and how it all fits together.
We chose Perl for  its strengths in text manipulation and its
interface with MySQL, our database software.

As explained in Section~\ref{S:math}, we need to generate a system of polynomials and compute
an eliminant, many times.
In previous computations, over $97\%$ of computer resources were spent on computing
eliminants. 
We need the computation to be efficient and to run on freely available software.
The methods of choice for elimination are algorithms based on Gr\"obner bases.
For this, we tested three Gr\"obner basis packages
(\cocoa~\cite{CoCoA}, Macaulay~2~\cite{M2}, and Singular~\cite{Singular}) on a
suite of representative problems.
When we made our choice of elimination software in the Autumn of 2007, Singular was by far the
fastest. 

Given an eliminant, we need to determine its number of real roots, quickly and reliably.
This requires a symbolic algorithm based on Sturm sequences~\cite{BPR}, and we needed
software that we could install on our many different computers.
While Maple is proprietary software, it has the fastest and most reliable
routine, {\tt realroot}, for counting real roots among the software we tested.
Maple was also installed on the computers we planned to use and we trust {\tt realroot}
completely, having used it on several billions of previous computations.

The mathematical routines of elimination and real root counting are symbolic (i.e., exact)
algorithms. 
We know of no satisfactory parallel implementations, so we achieve parallelism by running
different computations on different CPU cores.
\smallskip

When our software (a Perl script) is run, it queries the database for a Schubert
problem to work on and then writes a Singular input file to generate the 
desired polynomial systems and perform the eliminations.
As it writes this file, Perl selects random subsets $T$ of our master list of 111 rational
numbers, (re)orders them to make secant flags, and computes (and stores) the overlap number for
each polynomial system.
After the file is written, Perl calls Singular to run this file to 
compute the eliminants and write them to a file.
Perl then uses that output file to create an input file for Maple, which it 
calls Maple to run.
Maple determines the number of real roots of each eliminant, writing that to a file. 
Finally, Perl reads Maple's output, pairs the numbers of real roots with the corresponding
overlap number, posts these results to the database, and updates the state of the computation.

%
\subsection{Database}
A database is just a collection of {\sf tables} containing data, together with an interface
that allows efficient queries.
We designed a database to organize this computation.
It contains the Schubert problems to be studied, the results (e.g.~Table~\ref{T:Sample_table}),
and much else  in between.
At all times, the database contains a complete snapshot of the calculation.
Despite the size of this computation, the database is quite small, about 750 Kilobytes.
We briefly explain some of the more important tables in our database and their role in this
experiment.

{\sf Points} contains the master list of 111 numbers used to construct secant flags.
It is never altered.

{\sf SchubertProblems} contains the list of all Schubert problems we intend to study.
Section~\ref{S:load} explains how we add problems to the database.

{\sf Requests} keeps track of how much computation is to be done on each Schubert
problem and what has been started.
The administrators manually update {\sf Requests} to request more computation for particular
problems, and the Perl script updates {\sf Requests} when beginning a new computation on a
Schubert problem. 

{\sf Results} stores the frequency table of real solutions vs.~overlap number and the amount of
computing for each Schubert problem.
The Perl script updates {\sf Results} after successfully completing a run.
This table contains the information that our PHP web pages display.

{\sf RunningInstance} contains a list of the computations that have started but have yet to be
completed.  
We describe its function in Section~\ref{S:robustness}.

%
\subsection{Packets}\label{S:packets}
An important technical aspect of this computation is how we parcel out our
computations to individual computers.
There are many constraints.
Disk space is finite and large files are difficult to handle.
Some machines are available only for fixed time periods, and we must efficiently schedule their
use.
Networks and servers have fixed capacity, so database queries should be kept to a
minimum. 
Additionally, our computations require vastly different resources, with
some Schubert problems taking less than $0.033$ GigaHertz-seconds per instance while others
we studied require in excess of $40,000$ GigaHertz-seconds per instance.

To balance these constraints, we divide the computation of each Schubert problem into 
units that we call \DeCo{{\sl packets}}.
Each packet consists of between five and $50,000$ instances (one to $10,000$ choices of
the set $T$), and ideally requires about 1 hour of computation. 
The packets are processed through one or more Singular/Maple input files, none containing more
than 500 polynomial systems.
When a computer queries the database for a problem, it is given a single packet to
compute. 

The database stores the size and composition of the packets (which is set when the problems
are loaded into the database), and all information it records on the amount of computation
is denominated in these packets.

Packets for computationally-intensive Schubert problems require more than one hour
of computing.
Schubert problems are sorted by the expected time required for a packet, and this
is used in job scheduling to optimize performance.
The largest computations are performed on machines with no limit on their availability and the
others are parceled out according to the fit between the expected length of computation and the
computer's availability. 

%
\subsection{Loading Problems}\label{S:load}

Schubert problems are loaded into the database and the parameters of the computation are set
using different software than the main computation.
We have code to generate all Schubert problems on a given Grassmannian, determining the
number of solutions to each Schubert problem.
This uses a Gr\"obner bases computation in a standard
presentation of the cohomology ring, together with the Giambelli formula for Schubert
classes~\cite{Fu97}. 

A Perl script tests a subset of these problems to determine if the
computation is feasible.
An administrator selects feasible problems to
load into the database with a 
software tool that runs several instances of each problem and decides, based upon the
length of the computation, how to divide the computation into packets.
It writes these data into the database and records its work in a log file.

%
%
\section{Computing on a Beowulf cluster}\label{S:Calclabs}

A \DeCo{{\sl Beowulf cluster}} is a simple way to organize computers to work together
in which one machine (the server) communicates with the others (its clients), but there is no
communication between the clients. 
This model is optimal for performing many independent computations, for example, when running
several computations (e.g. computing Gr\"obner bases) in parallel.
It is a perfect match for our computational needs, and
most of our experiment runs on a Beowulf cluster.
We describe our cluster, its job scheduling, and how we
organized our use of this resource.

%
%
\subsection{The Calclabs cluster}

In Section~\ref{S:hardware} we mentioned our use of instructional computers at Texas A\&M
University which are collectively called the Calclabs~\cite{Calclabs}, as they are primarily
used in Calculus classes.
The Calclabs consist of 191 computers in five instructional labs and 12 in another lab.
Johnson installed the open source batch job scheduler Torque Resource Manager~\cite{Torque}
on these computers, which are the clients, and on a server.
Users log in to the server to submit jobs to a queue from which jobs are given to computers as
they becomes available.

Jobs are submitted to the queue with a specified time limit,
both for administration and because each computer is available only for a limited time period 
(typically nights, weekends, and holidays).
Jobs exceeding their time limit are terminated.
A computer is given the first job (if any) whose time limit does not exceed its availability. 
As described in Section~\ref{S:packets}, we sort Schubert problems by the expected
time to compute a packet to optimize this aspect of the scheduler.

%
%
\subsection{Scripting and cron}

While we monitor the progress of our computation on the Calclabs 
and sometimes submit jobs to the queue manually, we have largely automated the administration
of this computation with the Unix utility cron.
Cron executes scheduled tasks and is ideal for performing this administration.

We have set up cron on our account on the server to run a shell script which monitors the
queue, submitting jobs when the queue runs down.
It does this intelligently, ensuring that the queue contains packets of differing lengths,
tailored to the available computing.
This runs once per hour to keep the queue well-stocked.
Other administrative tasks (rotating logs, deleting old temporary files and archiving the
database) are performed once each day.

Since the scheduler submits one job to each machine, but our software runs on a single
core, the jobs are themselves shell scripts that run one copy of our software for each CPU core
on the given machine.

%
%
\section{Maintaining data quality}\label{S:quality}

An essential requirement in experimental science is that results are reproducible.
This is easy to ensure in computational experiments by using 
deterministic algorithms reliably implemented in software.
A second requirement is proper experimental design to ensure that a representative
sample has been tested.
Computational experiments may marry these two requirements by using (pseudo) random number
generators with fixed seeds for random sampling, and storing the seeds.

We explain our choices for experimental design and how we ensure the
reproducibility of our computation.
In principle, this experiment could be rerun, recreating every step in every computation.
This repeatability was essential for software development and testing, for checks on the
integrity of our data, and it will allow us to rerun a large part of the experiment using
different software on a different cluster.
With a computation of this complexity, failures of the software and networks are inevitable.
We explain how we recover from such failures, both those we anticipate and those that we
do not.

%
%
\subsection{Experimental design}
Schubert problems come in a countably infinite family with only a few tens of thousands small
enough to model on a computer.
We study many of the computable Schubert problems, and
for each, we test thousands to millions of instances of the Secant Conjecture.

Previous computations have shown the value of such indiscriminate testing.
The seminal example of the Monotone Conjecture (the cover illustration for the issue
of Experimental Mathematics in which the paper~\cite{RSSS} appeared) was tested 
late in that experiment, and only after the undergraduate member of that team asked why we were
omitting it.
(Sottile mistakenly thought it would be uninteresting.)
Also, extensive initial tests of the Shapiro Conjecture for flag manifolds appeared to affirm
its validity (it is in fact false).
Later was it realized that, by poor experimental design, only cases of 
the Monotone Conjecture had been tested, thereby overlooking counterexamples to the 
Shapiro Conjecture.

We kept these lessons in mind when designing the current experiment.
We were indiscriminate in selecting problems, studying all Schubert problems on
Grassmannians in 4-, 5-, and 6-dimensional space, as 
well as on $G(2,7)$ and $G(5,7)$, where $G(k,n)$ is the Grassmannian of $k$-planes in
$n$-space. 
We have also studied many computable problems on $G(3,7)$ and will study many on $G(4,7)$,
$G(5,7)$, $G(2,8)$, $G(3,8)$, and $G(4,8)$.
For these last six Grassmannians, we are computing a random selection of problems.
While it is hard to be precise, of the 7286 Schubert problems on $G(4,8)$, we estimate
that 3000 could be studied with our software.
The rest are too large to compute in a reasonable time or are infeasible.

For a problem involving Schubert conditions $(w_1,\dotsc,w_s)$, there are 
$s!$ ways to order the intervals for secants.
Our software randomly reorders the conditions before constructing secant
flags, to remove bias from the given ordering.
More serious is the question of how uniformly we are selecting from among all 
secant flags.
We do not have a satisfactory answer to this.
While one may believe that random subsets of our 111 master numbers (shown 
in Figure~\ref{F:RP1}) are fairly uniform modulo the action of $PSL(2,\R)$, we instead
offer experience gained in the previous experiment studying the Monotone Conjecture~\cite{RSSS}.
There, the results of a computation (e.g.~verifying the conjecture and an
inner border as in Table~\ref{T:Sample_table})
did not appear to depend upon how we selected subsets of a master set of numbers.
The selections included such schemes as all subsets of the numbers $\{1,\dotsc,10\}$, or
random subsets of the first 20 prime numbers, or random subsets of all rational numbers
$p/q$ where $(p,q)$ are the integer points closest to
$\left(101\cos(\frac{\pi}{40}n),101\sin(\frac{\pi}{40}n)\right)$, for $n=1,\dotsc,40$. 
(This last scheme is likely nearly uniform.)

%
%
\subsection{Reproducibility and seeds}
Random choices are made with the help of a pseudorandom
number generator.
Its output depends deterministically, but to all appearances unpredictably,
on a state variable called a \DeCo{{\sl seed}}, which is
deterministically updated after each call.
Thus two sequences of calls to a pseudorandom number generator beginning with different
seeds give unrelated sequences of integers, but if the seeds are 
the same, the sequences are identical.

We take advantage of this by generating an initial seed for each
Schubert problem  before computing its first packet.
This initial seed is determined by the current state of the computer.
When packet $n$ is begun, the seed is set to this initial seed and $n$
calls are made to the pseudorandom number generator to set
the seed for that packet.
In this way, the computation of a given Schubert problem is
completely determined given this initial seed, which is stored 
in our database.

%
%
\subsection{Independent Checks}
This exact reproducibility of results in a computational experiment is much stronger than
the notion of reproducibility for statistical results, and is a feature that we exploit.
We use it for software development to test upgrades
and to ensure that the software runs properly on different machines.
For this, we simply copy some problems and their initial seeds to an empty database, 
run all requested computations, and then compare the new results with the old results.
(They have always agreed.)

We have a software tool to automate this process and now use it on individual
laptops to rerun the computation for some Schubert problems to validate the data
for these problems from the initial run.
We understand that GIMPS also uses such double-checking for validation.

More interesting, and we believe unprecedented in mathematical experimentation, we are starting
to use this feature to rerun a large segment of the calculation on a different
cluster~\cite{Brazos} running a different version of Linux and different hardware and also
using different software for our basic mathematical routines. 
We use Macaulay~2 for elimination in place of Singular and SARAG~\cite{SARAG} in place of
Maple for counting real roots of univariate polynomials.
Besides providing an independent check on the data we generate, this will also give a direct
comparison of the efficacy of different implementations of these basic mathematical routines.

%
%
\subsection{Robustness}\label{S:robustness}
As with any complicated task, we cannot avoid the unexpected (the unknown unknowns), and 
have designed our software to recover from the many different failures that inevitably
occur.
For this, we have several interlocking systems to prevent corrupted calculations from
being entered in the database.
We also rerun corrupted packets, and even rerun all or part of a
Schubert problem whose data appear suspect. 

First, our software has checks to ensure that tasks (connecting with the database,
calling client programs, and reading/writing data in files) are
successfully executed, and which terminate its running when something untoward is detected.
Other errors, such a network errors or unmounted file systems, cause noisier failures, which
are captured by our log files for possible diagnosis.
There are even less graceful failures in a computation, such as power outage, termination of jobs by
the job server, or simple human error.
All of these abort the computation of a packet and therefore lead to packets
whose computation has started, but whose results have not been submitted to the database.

The table {\sf RunningInstance} in our database keeps track of packets whose computation has
started but not finished, together with the expected completion time.
Packets that have been terminated in any way are recognized by
having an expected completion time that has long passed.
When our software queries the database for a problem to work on, it first checks for any such 
overdue packets.
If one is found, it deletes that record from {\sf RunningInstance} and creates a new record
corresponding to this new computation.
Otherwise it finds a fresh packet to compute, creating a record in {\sf RunningInstance}.
Upon successful completion, these records are deleted.
One possible graceful failure is for a computation to end, but discover that its
record in {\sf RunningInstance} has been deleted and superseded by another---preventing a
second submission of the same data to the database.

While this method works for the few to hundreds (out of thousands) of packets each day that fail
before successful completion, sometimes our data becomes corrupted, or possibly corrupted.
We also have a software tool that finds the most recent database backup where that Schubert
problem is uncorrupted and restores that Schubert problem to this previous state.
Thus we simply recompute all or part of the computations for that Schubert problem.

%
%
\section{Conclusions and future work}\label{S:conclusions}

We plan to continue this work of mathematical discovery through
advanced computing.
At the conclusion of this experiment in December 2009, we will write a full paper describing its
mathematical results.
As of May 2009, the Secant Conjecture was verified in each of the over 250 million
instances we checked.
In 2010 we plan to start a related experiment, testing a
common generalization of the Monotone and Secant Conjectures.
This will last about one year.
While there is much more to be discovered studying these variants of the Shapiro Conjecture, 
we plan a long-term, multifaceted, and systematic study of Galois groups of Schubert
problems, building on the work in~\cite{BV,LS_Galois}.

We mention one side benefit from this computation.
In March 2009, after sharing timing data from a benchmark computation with Mike Stillman, a
developer of Macaulay~2, he rewrote some code that improved its running by several
orders of magnitude.\smallskip 

We have described how and why we set up, organized, and are running a very large computational
experiment to study a conjecture in pure mathematics, and how it is possible to 
harness underused yet widely available computing resources for mathematical research.
We believe this model---a team-based approach to designing and monitoring a large
computational experiment---can fruitfully be replicated in other settings and we encourage you
to try it.

\providecommand{\bysame}{\leavevmode\hbox to3em{\hrulefill}\thinspace}
\providecommand{\MR}{\relax\ifhmode\unskip\space\fi MR }
\providecommand{\MRhref}[2]{%
  \href{http://www.ams.org/mathscinet-getitem?mr=#1}{#2}
}
\providecommand{\href}[2]{#2}


\end{document}